\documentclass{scrartcl}
\usepackage{authblk}
\usepackage{amsmath}
\usepackage{amssymb}
\usepackage{mathrsfs}
\usepackage[T1]{fontenc}
\usepackage{pxfonts}

\usepackage{comment}
\usepackage{color}
\usepackage{bm}
\usepackage{float}
\usepackage{graphicx}
\usepackage{wrapfig}
\usepackage{url}

\newcommand{\argmax}{\mathop{\rm argmax}\limits}

\title{An implicit method for the finite time horizon Hamilton-Jacobi-Bellman quasi-variational inequalities}
\date{\today}
\author[$\dag$]{Masashi Ieda\thanks{ieda@craft.titech.ac.jp}}
\affil[$\dag$]{
	Graduate School of Innovation Management \authorcr
	Tokyo Institute of Technology \authorcr
	2-12-1 Ookayama, Meguro-ku, Tokyo, Japan
}

\date{}

\begin{document}

\maketitle

\begin{abstract}
We propose a new numerical method for solving the Hamilton-Jacobi-Bellman quasi-variational inequality
 associated with the combined impulse and stochastic optimal control problem over a finite time horizon.
Our method corresponds to an implicit method in the field of numerical methods for partial differential equations,
 and thus it is advantageous in the sense that the stability condition is independent of the discretization parameters.
We apply our method to the finite time horizon optimal forest harvesting problem, which considers exiting from the forestry business at a finite time.
We show that the behavior of the obtained optimal harvesting strategy of the extended problem coincides with our intuition.
\end{abstract}

\section{Introduction}
Solving the Hamilton-Jacobi-Bellman quasi-variational inequality (HJBQVI) is one of the most challenging issues
 in the stochastic optimal control problem.
The HJBQVI is associated with the combined impulse and stochastic optimal control, which can be used to formulate
 a system which changes drastically under our control.
The combined stochastic optimal control is a widely applicable framework.
Some of the literature that deals with applications to mathematical finance is as follows.
 Pliska and Suzuki \cite{Pliska2004}, Palczewski and Zabczyk \cite{Palczewski2005}
  and Kharroubi and Pham \cite{Kharroubi2010a} treat portfolio optimization with transaction costs;
 Mundaca and {\O}ksendal \cite{Mundaca1998} and Cadenillas and Zapatero \cite{Cadenillas2000} study control of the exchange rate by the Central bank;
 Korn \cite{Korn1999} provides an overview relating to applications of the impulse control.
Applications to other areas, such as electricity management, problems of maintenance and quality control and information technology,
 are found in Bensoussan and Lions \cite{bensoussan1984impulse} and the references therein.

A common approach to applying the impulse control framework has the issue that problems or models are formulated to have an analytical solution.
For the HJBQVI associated with the one-dimensional infinite horizon combined stochastic optimal control problem,
 the smooth-fit technique is an established method for obtaining the solution.
However, this technique is not valid for the general HJBQVI, and to the best of our knowledge there is no established method for the general HJBQVI.
Hence, developing a numerical method for solving the HJBQVI is an important avenue of research. 

The numerical approaches to solving the HJBQVI are categorized into two types according to the time horizon, which is either infinite or finite.
We first discuss the infinite time horizon.
Bensoussan and Lions \cite{bensoussan1984impulse} approximate the impulse control problem using iterations of the optimal stopping problem,
 and hence the HJBQVI is translated to HJB variational inequalities (HJBVIs).
The numerical method for the HJBVI is well-studied due to the motivation
 of pricing American options in mathematical finance.
An alternative method is proposed in Chancelier et. al. \cite{Chancelier2006}. 
In this paper, the authors provide a solution method for a fixed-point problem consisting of a contractive operator
 and a non-expansive operator, and the numerical algorithm for the infinite time horizon HJBQVI appears as an application.
They discretize the HJBQVI using a finite difference scheme and the discretized HJBQVI is converted to an equivalent fixed-point problem that is solvable by their algorithm.

In the case of the finite time horizon, we can use backward induction
 in a fashion similar to that in other optimal control problems with a terminal condition.
For instance, Chen and Forsyth \cite{Chen2008} solve the finite time horizon HJBQVI associated with
 the annuity pricing problem with a guaranteed minimum withdrawal benefit.
They construct a discretized equation that is consistent with the original HJBQVI
 under the constraint that the time grid size is small enough to satisfy a certain condition depending on the space grid size.
We note that this method is dependent on the model, and thus it appears to be difficult to apply directly to other problems.

In the present paper, we propose a new numerical algorithm for solving the finite time horizon HJBQVI.
Compared with the infinite time horizon case, we must treat the time variable carefully.
We provide an appropriate transformation for the HJBQVI considered,
 and develop a numerical scheme that does not require a relationship between the time and special discretizations.
The key points of the present study are:
 (i) we discretize the HJBQVI with the forward difference for the time grid; and
 (ii) the impulse control part is kept in the same time step.
Hence the present method is an implicit method.

The main advantage of the present work is providing a numerical scheme independent of the specific model formulation.
The stability of the numerical scheme is especially sensitive in the finite time horizon case.
Previous research such as \cite{Chen2008} proves the stability for each of the concerned problems.
We give a proof of stability in the general framework, and hence
 our method can be applied without a stability proof.
To support the advantages of our method from the perspective of implementation,
 we provide a detailed procedure of the algorithm in matrix form.

This paper is organized as follows.
We present the mathematical formulation of the HJBQVI in Section \ref{sec:formulation}.
The goal of this section is to display the discretized HJBQVI in matrix form.
Section 3 provides details of our algorithm.
In Section 3.1, the matrix form of the HJBQVI obtained in Section \ref{sec:formulation} is translated into an equivalent fixed-point problem.
We describe the procedure for solving the fixed-point problem in detail from the view point of the computational implementation in Section \ref{sec:algo}.
In Section \ref{sec:num_res}, we apply the proposed method to the optimal forest harvesting problem.
This problem determines the optimal harvesting strategy for ongoing forestry.
We use the mathematical formulation of the problem proposed by Willassen \cite{Willassen1998}
 based on an infinite time horizon impulse control framework 
 where analytical solutions of the value function and optimal strategy are provided.
We introduce a terminal time representing the time of exit from the forestry business,
 which turns the above problem into a finite time horizon problem. 
In this case, the analytical solution is unavailable and thus we solve it numerically using the proposed algorithm.

\section{Mathematical formulation}
\label{sec:formulation}
\subsection{Model}
We consider the following combined stochastic and impulse control problems over a finite time horizon $[0,T]$ with the performance criterion
\begin{equation}
\label{eq;performance_cri}
J^w_t(x) = \mathbb{E}\left[ \left. \int_t^T f(s,X_s^{w},u_s) dt +g(X_T^{w}) + \sum_{t<\tau_j<T} K(\tau_j^-, X_{\tau_j^-}^{w},\zeta_j ) \right| X^w_t = x \right].
\end{equation}
The controlled process $\left\{ X_t^{w} \right\}_{t\geq 0}$ is governed by
\begin{align}
	&\begin{cases}
		\displaystyle dX_t^{w} = \mu ( t, X_t^{w},u_t ) dt + \sigma ( t, X_t^{w},u_t ) dW_t, \quad \tau_j \leq t < \tau_{j+1},\\
		X_{\tau_{j+1}} = \Gamma( X_{\tau_{j+1}^-}, \zeta_{j+1} ), \quad j=0,1,2,\cdots,
	\end{cases}
	 \label{eq:def_ctrl_process}
\end{align}
where $W_t$ is a $d$-dimensional Brownian Motion on a filtered probability space $( \Omega, \mathcal{F}, \left\{ \mathcal{F}_t \right\}_{t\geq0}, \mathbb{P})$,
 $\mu:[0,T]\times\mathbb{R}^n\times \mathbb{U} \rightarrow \mathbb{R}^n$, 
 $\sigma:[0,T]\times\mathbb{R}^n\times \mathbb{U} \rightarrow \mathbb{R}^{n\times d}$,
 $\Gamma : \mathbb{R}^n\times\mathcal{Z}\rightarrow\mathbb{R}^n$,
 $\mathbb{U}\subset \mathbb{R}^l$,
 and  $\mathcal{Z}\subset \mathbb{R}^l$.
The combined control $w$ consists of the Markov control strategy $u=\{u_t\}_{t\geq0}$, 
 a $\mathbb{U}$-valued stochastic process of the form $u_t=\alpha(t,X^w_t)$ for some function $\alpha:[0,T)\times\mathbb{R}^n\rightarrow \mathbb{U}$,
 and the impulse control strategy $v=\{(\tau_j,\zeta_j)\}_{j=1}^\infty$. 
Here $\tau_0=0$, $\tau_1<\tau_2<\cdots$ are $\mathcal{F}_t$-stopping times
 and $\zeta_j\in\mathcal{Z}$, $j\geq 1$, are $\mathcal{F}_{\tau_j}$-measurable random variables.
The performance is measured using the following three functions:
 the profit rate function $f:[0,T]\times\mathbb{R}^n\times \mathbb{U} \rightarrow\mathbb{R}$ ,
 the bequest function $g:\mathbb{R}^n\rightarrow\mathbb{R}$,
 and the intervention profit function $K:[0,T),\mathbb{R}^n\times\mathcal{Z}\rightarrow\mathbb{R}$.

We denote by $\mathcal{W}$ the set of admissible combined controls, that is,
$w\in\mathcal{W}$ satisfies: (i) a unique strong solution of the SDE (\ref{eq:def_ctrl_process}) with control $w$ exists;
 (ii) $\lim_{j\rightarrow\infty} \tau_j=T$ almost surely.
We also assume that for $w\in\mathcal{W}$,
$$
	\mathbb{E} \left[ \int_0^T \left| f(t,X_t^{w},u_t) \right| dt \right] < \infty,\;
	\mathbb{E} \left[ \left| g(X_T^{w}) \right| \right] < \infty,\;
	\mathbb{E} \left[ \sum_{\tau_j<T} \left| K(\tau_j^-, X_{\tau_j^-}^{w},\zeta_j ) \right| \right] < \infty.
$$

The value function corresponding to problem (\ref{eq;performance_cri}) is defined by $V_t(x)=\sup_{w\in \mathcal{W}}J^w_t(x)$, $(t,x)\in[0,T]\times\mathbb{R}^n$.
Hence the corresponding HJBQVI is given by (see, for example, \cite{oksendal2007applied})
\begin{equation}
\max\left(
	\sup_{\alpha\in\mathbb{U}} \left\{ \partial_t V_t(x)+\mathscr{L}_t^{\alpha} V_t(x) + f(t,x,\alpha) \right\},
	\sup_{\zeta\in\mathcal{Z}} \left\{ V_t(\Gamma(x,\zeta))+K(t,x,\zeta) - V_t(x) \right\}
 \right)=0,
\label{eq:general_qvi}
\end{equation}
with terminal condition $V_T(x)=g(x)$, $x\in\mathbb{R}^n$, where $\partial_t$ is the partial differential operator with respect to $t$, and
$\mathscr{L}_t^{\alpha}$ is the infinitesimal generator of the process $X^w$ at time $t$:
$$
\mathscr{L}_t^{\alpha} \Psi(x) = \mu(t,x,\alpha)^\top \partial_x \Psi(x)
 + \frac{1}{2} \mathrm{Tr} \left( (\sigma\sigma^\top)(t,x,\alpha)\partial_x^2\Psi(x) \right)
$$
for $\Psi\in C^{2}(\mathbb{R}^n)$.
Here $\partial_x^j$ is the $j$-th order partial differential operator with respect to $x$
 and $\top$ indicates transposition.
The continuation set $\mathcal{D}_t$ at time $t$ is defined by
$$
\mathcal{D}_t = \left\{ x\in \mathbb{R}^n \left| \sup_{\zeta\in\mathcal{Z}} \left\{ V_t(\Gamma(x,\zeta))+K(t,x,\zeta) - V_t(x) \right\} <0  \right. \right\}.
$$
We impose the condition $\mathcal{D}_t \neq \emptyset$, $t\in[0,T)$ on the control set, that is, a control which intervenes in all regions is not admissible.
A simple way to satisfy this condition is that the intervention times $\{\tau_j\}_{j=1}^\infty$ satisfy the following statement:
for each $t \in (0,T)$, there exists $x^R \in \mathcal{S}$ such that $X_t=x^R$ and $\tau_{j+1}\neq t$
 where $j = \argmax_{i} \left\{ \tau_i <  t \right\}$.
We are usually able to find such an $x^R$ from the problem formulation: an example of this is given in Section \ref{sec:ofh_form}.


As with other numerical problems, we face the issue of the computer not dealing well with the unbounded domain.
Since our state process $X_t^{w}$ leaves any bounded region until the termination time $T$ with non-zero probability, this problem is unavoidable.
To cope with it, we restrict ourselves to problems with a bounded domain.
Let $\mathcal{S}\subset\mathbb{R}^n$ be the domain and $\partial\mathcal{S}$ be the boundary.
We assume that (i) $\mathcal{D}_t\cap\mathcal{S} \neq \emptyset$, $t\in[0,T)$;
(ii) the intervention function $\Gamma$ satisfies $\Gamma:\mathcal{S}\rightarrow\mathcal{S}$.
The boundary condition is given by the function $\psi:[0,T)\times\partial\mathcal{S}\rightarrow\mathbb{R}$.

\subsection{Discretization}
We discretize the QVI (\ref{eq:general_qvi}) using the standard finite difference scheme with the central difference.
Let $\delta_t,\delta = (\delta_1,\cdots,\delta_n)^\top$ be the finite difference steps with respect to $t$ and $x$.
We denote by $\mathcal{S}_\delta$ the spacial grid, so that $\mathcal{S}_\delta=\mathcal{S}\cap\prod_{i=1}^n(\delta_i\mathbb{Z})$.
The discretized boundary $\partial\mathcal{S}_\delta$ may also be represented by $\partial\mathcal{S}_\delta=\partial\mathcal{S}\cap\prod_{i=1}^n(\delta_i\mathbb{Z})$.
For convenience we introduce symbols for the time grid points, $\{t_i\}_{0\leq i \leq N^t}$, $t_i\in[0,T]$,
 and for the spacial grid points including the boundary, $\{x_i\}_{1 \leq i \leq N^x+\bar{N}^x}$, $x_i\in\mathcal{S}_\delta\cup\partial\mathcal{S}_\delta$,
where $N^t$ is the number of time grid points, $N^x$ is the number of spacial grid points and
 $\bar{N}^x$ is the number of discretized boundary points.
We note that the subsequences $\{x_i\}_{1 \leq i \leq N^x}$ and $\{x_i\}_{N^x < i \leq N^x+\bar{N}^x}$
 represent the spacial (internal) grid points and discretized boundary points respectively.
Furthermore, the time grid points are defined sequentially as $t_0=0$, $t_1=\delta_t$, $\cdots$, $t_{N^t}=N^t\delta_t=T$.
We implicitly assume that $\delta_t$ is the number supporting the existence of $N^x\in\mathbb{N}$ such that $N^t\delta_t=T$.
Let $\mathcal{Z}_\delta$ be a discretized set of the impulse control, and 
assume that:
(i) $\Gamma:\mathcal{S}_\delta \times \mathcal{Z}_{\delta} \rightarrow \mathcal{S}_\delta$,
(ii) we can take space indices which allow to define an integer function
 $\eta:\{1,\cdots,N^x\}\times\mathcal{Z}_\delta\rightarrow\{1,\cdots,N^x\}$ such that
\begin{equation}
\Gamma(x_i,\zeta) =x_{\eta(i,\zeta)}, \quad \eta(i,\zeta)<i, \quad \mathrm{with} \quad i\in\{1,\cdots,N^x\},\; \zeta\in\mathcal{Z}_\delta,
\label{eq:ass_gamma}
\end{equation}
and $x_1 \in \mathcal{D}_t$.

The discretized QVI of the QVI (\ref{eq:general_qvi}) is defined as follows:
\begin{align}
	&\begin{cases}
		\displaystyle\max\left(
			\sup_{\alpha \in\mathbb{U}} \left\{ \frac{\Phi^{k+1}(x)-\Phi^k(x)}{\delta_t}+\bar{\mathscr{L}}^{\alpha,k} \Phi^k(x) + f(t_k,x,\alpha) \right\},\right. \\
		\displaystyle \qquad\qquad\qquad \left. \sup_{\zeta\in\mathcal{Z}_\delta} \left\{ \Phi^k(\Gamma(x,\zeta)) +K(t_k,x,\zeta) -\Phi^k(x) \right\}
 		\right)=0,\quad x\in\mathcal{S}_\delta,\\
 		\displaystyle\Phi^k(x) = \psi(t_k,x),\quad x\in \partial\mathcal{S}_\delta,
 	\end{cases}
\label{eq:discretized_general_qvi}
\end{align}
for $k\in\{0,\cdots,N^t-1\}$,
with terminal condition $\Phi^{N^t}(x)=g(x)$, $x\in\mathcal{S}_\delta\cup\partial\mathcal{S}_\delta$,
where $\Phi^k:\mathcal{S}_\delta\cup\partial\mathcal{S}_\delta\rightarrow\mathbb{R}$
 and $\bar{\mathscr{L}}^{\alpha,k}$ is the operator such that 
\begin{align*}
	\bar{\mathscr{L}}^{\alpha,k} \Psi(x) =&
		\displaystyle\Psi(x)\left\{ \sum_{i=1}^n \frac{-(\sigma\sigma^\top)_{ii}}{\delta_i^2} + \sum_{j \in \mathcal{J}(i)}\frac{|(\sigma\sigma^\top)_{ij}|}{2\delta_i\delta_j} \right\}\\
 		&\displaystyle \; + \frac{1}{2} \sum_{i=1}^n \sum_{\kappa=\pm1} \Psi(x+\kappa\delta_ie_i) \left\{ \frac{-(\sigma\sigma^\top)_{ii}}{\delta_i^2}
 				 - \sum_{j \in \mathcal{J}(i)}\frac{|(\sigma\sigma^\top)_{ij}|}{\delta_i\delta_j} + \kappa\frac{\mu_i}{\delta_i} \right\}\\
 		&\displaystyle \; + \frac{1}{2} \sum_{i=1} \sum_{j \in \mathcal{J}(i)} \sum_{\kappa,\lambda=\pm1} \Psi(x+\kappa\delta_ie_i+\lambda\delta_je_j)
 				 \frac{(\sigma\sigma^\top)_{ij}^{[\kappa\lambda]}}{\delta_i \delta_j}, \quad x \in \mathcal{S}_\delta,
\end{align*}
for $\Psi:\mathcal{S}_\delta\cup\partial\mathcal{S}_\delta\rightarrow\mathbb{R}$.
Here we have used the notation
\begin{align*}
	(\sigma\sigma^\top)_{ij}^{[\kappa\lambda]} =
	&\begin{cases}
		\displaystyle \max \left( 0, (\sigma\sigma^\top)_{ij} \right), & \kappa\lambda=1,\\
		\displaystyle - \min \left( 0, (\sigma\sigma^\top)_{ij} \right), & \kappa\lambda=-1,
	\end{cases}\\
	\mathcal{J}(i) =& \{1,\cdots,n\}\setminus\{i\},
\end{align*}
and have omitted the arguments of $(\sigma\sigma^\top)$ and $\mu$, that is, $(\sigma\sigma^\top)=(\sigma\sigma^\top)(t_k,x,\alpha)$ and $\mu=\mu(t_k,x,\alpha)$.

We use the central difference to obtain the discretized QVI (\ref{eq:discretized_general_qvi}).
If we use the one-sided difference, the discretized operator $\bar{\mathscr{L}}^{\alpha,k}$ becomes
\begin{align*}
	\bar{\mathscr{L}}^{\alpha,k} \Psi(x) =&
		\displaystyle\Psi(x)\left\{ \sum_{i=1}^n \frac{-(\sigma\sigma^\top)_{ii}}{\delta_i^2} 
		+ \sum_{j \in \mathcal{J}(i)}\frac{|(\sigma\sigma^\top)_{ij}|}{2\delta_i\delta_j}
		- \frac{|\mu_i|}{\delta_i}		
		 \right\}\\
 		&\displaystyle \; + \frac{1}{2} \sum_{i=1}^n \sum_{\kappa=\pm1} \Psi(x+\kappa\delta_ie_i)
 		 \left\{ \frac{-(\sigma\sigma^\top)_{ii}}{\delta_i^2}
 				 - \sum_{j \in \mathcal{J}(i)}\frac{|(\sigma\sigma^\top)_{ij}|}{\delta_i\delta_j} 
 				 + \frac{\mu_i^{[\kappa]}}{\delta_i}
 		\right\}\\
 		&\displaystyle \; + \frac{1}{2} \sum_{i=1} \sum_{j \in \mathcal{J}(i)} \sum_{\kappa,\lambda=\pm1} \Psi(x+\kappa\delta_ie_i+\lambda\delta_je_j)
 				 \frac{(\sigma\sigma^\top)_{ij}^{[\kappa\lambda]}}{\delta_i \delta_j}, \quad x \in \mathcal{S}_\delta.
\end{align*}

Finally, we represent the discretized QVI using the central difference (\ref{eq:discretized_general_qvi}) in matrix QVI form:
\begin{align}
&\begin{cases}
\max\left( 
	\displaystyle\sup_{a\in\mathbb{U}^{N^x}} \left\{ \frac{\phi^{k+1}_i -\phi^k_i}{\delta_t} + \left( L^{a,k}\phi^k \right)_i  + f^{a,k}_i  \right\} \right.\\
	\displaystyle \left. \hspace*{2cm} ,\sup_{z\in\mathcal{Z}^{N^x}_\delta} \left\{ \left( M^{z}\phi^k \right)_i +K^{z,k}_i -\phi^k_i \right\}
\right)=0, \quad 1 \leq i \leq N^x,\\
\phi^k_i = \psi(t_k,x_i), \quad N^x < i \leq N^x+\bar{N}^x,
\end{cases}
\label{eq:matrix_qvi}
\end{align}
for $k\in \{0,\cdots, N^t-1\}$ with terminal condition $\phi^{N^t}_i = g(x_i)$, $i\in\{1,\cdots, N^x\}$,
where $\phi^k$ is an $(N^x+\bar{N}^x)$-dimensional vector such that $\phi^k_{i}=\Phi^k(x_i)$,
$a=\{ a_i \}_{1\leq i \leq N^x}$, $a_i\in\mathbb{U}$,
$z=\{ z_i \}_{1\leq i \leq N^x}$, $z_i\in\mathcal{Z}_\delta$,
$L^{\alpha_k,k}$ is an $N^x \times (N^x+\bar{N}^x)$ matrix such that, for $i\in\{1,\cdots,N^x\}$,
\begin{align*}
	L^{a,k}_{ij} =
	\begin{cases}
		\displaystyle \sum_{i^\prime=1}^n \frac{-(\sigma\sigma^\top)_{i^\prime i^\prime}  }{\delta_{i^\prime}^2}
		 	+ \sum_{j^\prime \in \mathcal{J}(i^\prime)}\frac{|(\sigma\sigma^\top)_{i^\prime j^\prime} |}{2\delta_{i^\prime}\delta_{j^\prime}},
			\hspace*{2.5cm} \textrm{if} \;  i=j,\\
		\displaystyle \frac{1}{2} \left\{ \frac{-(\sigma\sigma^\top)_{i^\prime i^\prime}  }{\delta_{i^\prime}^2}
 				 - \sum_{j^\prime \in \mathcal{J}(i^\prime)}\frac{|(\sigma\sigma^\top)_{i^\prime j^\prime}|}{\delta_{i^\prime} \delta_{j^\prime}}
 				 + \kappa\frac{\mu_{i^\prime}}{\delta_{i^\prime}}  \right\},
			\hspace*{0.1cm} \begin{array}{l}
					\textrm{if} \;   x_j= x_i+\kappa\delta_{i^\prime} e_{i^\prime}, \\
					\; \kappa\pm 1, \; i^\prime \in \{ 1,\cdots,n \}
				\end{array}\\
		\displaystyle \frac{1}{2}\frac{(\sigma\sigma^\top)_{i^\prime j^\prime}^{[\kappa\lambda]} }{\delta_{i^\prime} \delta_{j^\prime}},
			\hspace*{6cm} \begin{array}{l}
					\textrm{if} \;  x_j =x_i+\kappa\delta_{i^\prime} e_{i^\prime} + \lambda \delta_{j^\prime} e_{j^\prime},\\
					 \; \kappa,\lambda=\pm1,\\
					 \; i^\prime \in \{ 1,\cdots,n \},\\
					 \; j^\prime \in \mathcal{J}(i^\prime),
				\end{array}\\
		\displaystyle 0, \hspace*{8.5cm} \textrm{otherwise},
	\end{cases}
\end{align*}
where $f^{a,k}$ is an $N^x$-dimensional vector such that $f^{a,k}_i=f(t_k,x_i,a_i)$,
$M^{z}$ is an $N^x \times (N^x+\bar{N}^x)$ matrix such that $M^{z}_{ij} = \mathbf{1}_{\{ j = \eta(i,z_i) \}}$
and $K^{z,k}$ is an $N^x$-dimensional vector such that $K^{z,k}_i=K(t_k,x_i,z_i)$.
We have again omitted the arguments of $(\sigma\sigma^*)$ and $\mu$,
 that is, $(\sigma\sigma^\top)=(\sigma\sigma^\top)(t_k,x_i,a_i)$ and $\mu=\mu(t_k,x_i,a_i)$.

We assume that the functions $\mu$ and $\sigma$ satisfy the condition
\begin{align}
 |\mu_i(t,x,\alpha)| \leq \frac{(\sigma\sigma^\top)_{ii}(t,x,\alpha)}{\delta_i} - \sum_{j \in \mathcal{J}(i)}\frac{|(\sigma\sigma^\top)_{ij}(t,x,\alpha)|}{\delta_j},
 	\quad (t,x,\alpha) \in [0,T)\times\mathcal{S}_\delta\times\mathbb{U},
 \label{eq:ass_central}
\end{align}
so that the equation (\ref{eq:matrix_qvi}) will be stable (see Appendix \ref{sec:stability}).
We remark that if we use the one-sided difference, the condition (\ref{eq:ass_central}) is milder:
\[
 	0 \leq \frac{(\sigma\sigma^\top)_{ii}(t,x,\alpha)}{\delta_i} - \sum_{j \in \mathcal{J}(i)}\frac{|(\sigma\sigma^\top)_{ij}(t,x,\alpha)|}{\delta_j},
 		\quad (t,x,\alpha) \in [0,T)\times\mathcal{S}_\delta\times\mathbb{U},
\]
even though the convergence speed is slower.

\section{Algorithm}

\subsection{The equivalent fixed-point problem}
We convert the matrix QVI (\ref{eq:matrix_qvi}) to the following equivalent fixed-point problem:
\begin{align}
&\begin{cases}
\phi_i = \max\left( 
	\displaystyle\sup_{a\in\mathbb{U}^{N^x}} \left\{ \left( \bar{L}^{a,k}\phi^k \right)_i  + \bar{f}^{a,k}_i  \right\}
	 ,\sup_{z\in\mathcal{Z}^{N^x}} \left\{ \left( M^{z}\phi^k \right)_i +K^{z,k}_i \right\}
\right), \quad 1 \leq i \leq N^x,\\
\phi^k_i = \psi(t_k,x_i), \quad N^x < i \leq N^x+\bar{N}^x,
\end{cases}
\label{eq:conv_matrix_qvi}
\end{align}
where $\bar{f}^{a,k}$ is an $N^x$-dimensional vector and $\bar{L}^{a,k}$ is an $N^x \times (N^x+\bar{N}^x) $ matrix
defined respectively as follows:
\begin{align}
	&\begin{cases}
	 	\displaystyle \bar{f}^{a,k}_i = \frac{h \delta_t}{h+\delta_t} f^{a,k}_i + \frac{h}{h+\delta_t}\phi^{k+1}_i,\\
	 	\displaystyle \bar{L}^{a,k} = \frac{\delta_t}{h+\delta_t}(I + h L^{a,k}).
	\end{cases}
	\label{eq:bar_lf}
\end{align}
Here, $h$ is a positive number such that
\[
h \leq \inf_{a\in\mathbb{U}^{N^x}} \min_{1\leq i \leq N^x , 1 \leq k < N^t} \frac{1}{ \left| L^{a,k}_{ii} \right|},
\]
and $I$ is an $N^x \times (N^x+\bar{N}^x) $ matrix such that $I_{ij}=\delta_{ij}$ where $\delta_{ij}$ is the Kronecker delta.
Then the fixed-point problem (\ref{eq:conv_matrix_qvi}) satisfies the conditions for applying the method proposed by Chancelier et al. \cite{Chancelier2006}.
The details are discussed in Appendix \ref{sec:conv_mat_fix}.

\subsection{Procedures}
\label{sec:algo}
We are able to solve problem (\ref{eq:conv_matrix_qvi}) by the backward method.
In the following steps, we assume that we have already found $\phi^{k+1}$.
\begin{description}
\item[step1] set $\phi^k=\phi^{k+1}$,
\item[step2] search for the controls $\hat{a}^k$ and $\hat{z}^k$ such that
	$$
		\hat{a}^k = \argmax_{a\in\mathbb{U}^{N^x}} \left\{ \bar{L}^{a,k}\phi^k +\bar{f}^{a,k}\right\}, \quad
		\hat{z}^k = \argmax_{z\in\mathcal{Z}^{N^x}} \left\{ M^z\phi^k + K^{z,k} \right\},
	$$
and define an index set $\mathcal{I}^k$ such that
	$$
	\mathcal{I}^k = \left\{ i \in \{1,\cdots,N^x\} \left| (\bar{L}^{\hat{a}^k,k}\phi^k)_i +\bar{f}^{\hat{a}^k,k}_i \geq (M^{\hat{z}^k}\phi^k)_i + K^{\hat{z}^k,k}_i \right. \right\},
	$$
\item[step3] determine an $(N^x+\bar{N}^x)\times (N^x+\bar{N}^x)$ matrix $A$ and an $(N^x+\bar{N}^x)$-dimensional vector $b$ as follows:
	\begin{align*}
		&A_{ij} = \begin{cases}
			\bar{L}^{\hat{a}^k,k}_{ij} & \mathrm{if}\; i \in \mathcal{I},\\
			M^{\hat{z}^k}_{ij} & \mathrm{if}\; i \in \{1,\cdots,N^x\}\backslash \mathcal{I},\\
			0 & \mathrm{if} \; i \in \{N^x+1,\cdots,N^x+\bar{N^x}\},
		\end{cases}\\
		&b_i = \begin{cases}
			\bar{f}^{\hat{a}^k,k}_i & \mathrm{if}\; i \in \mathcal{I},\\
			K^{\hat{z}^k,k}_i & \mathrm{if}\; i \in \{1,\cdots,N^x\}\backslash \mathcal{I},\\
			\psi(t_k,x_i) & \mathrm{if} \; i \in \{N^x+1,\cdots,N^x+\bar{N^x}\},
		\end{cases}
	\end{align*}
\item[step4] solve the linear equation $(I-A)\phi^\prime=b$, where $\phi^\prime$ is an $(N^x+\bar{N}^x)$-dimensional vector
 and $I$ is an $(N^x+\bar{N}^x)\times (N^x+\bar{N}^x)$ identity matrix,
\item[step5] if $\max|\phi^\prime_i - \phi_i^k|$ exceeds the admissible error, replace $\phi^k$ by $\phi^\prime$ and go back to step2,
 else replace $\phi^k$ by $\phi^\prime$ and go to step6,
\item[step6] if $k\neq 1$, replace $k$ by $k-1$ and go back to step1, 
 else determine the Markov control $\hat{\alpha}$, the set $\hat{\mathcal{D}}_{t_k}$ and the impulse control $\hat{v}=(\hat{\tau_i},\hat{\zeta}_i)_{i\geq1}$ as follows:
	\begin{align*}
  		&\hat{\alpha}(t_k,x_i) = \hat{a}^k_i, \quad
  		\hat{\mathcal{D}}_{t_k}=\left\{ x_i \left| i \in \mathcal{I}^k \right. \right\},\\
		&\hat{\tau}_i = \min \left\{ t_k \in \{t_1,\cdots,t_{N^t} \}  \left| t_k > \hat{\tau}_{i-1}; X^{\hat{w}}_{t_k} \notin \hat{\mathcal{D}}_{t_k} \right. \right\},\\
		&\hat{\zeta}_i(x_j) = \hat{z}^{k^\prime}_j, \quad k^\prime  \in \{1,\cdots,N^t\} \;\mathrm{s.t.}\; {t_{k^\prime}=\hat{\tau}_i},
	\end{align*}
where $\hat{\tau}_0=0$ and $\hat{w}=(\hat{\alpha},\hat{v})$.
\end{description}
The control $\hat{w}$ obtained by the above procedure satisfies $J^{\hat{w}}_{t_k}(x_i) \geq J^w_{t_k}(x_i)$, $t_k \in \{t_0,\cdots,t_{N^t-1} \}$, $x_i\in\mathcal{S}_\delta$, $w\in\mathcal{W}$,
and hence $\hat{w}$ is an optimal control.

\section{Numerical results}
\label{sec:num_res}
We apply our method to the finite horizon optimal forest harvesting problem based on Willassen's formulation.
The mathematical formulations of the original and extended problems are described in Section \ref{sec:ofh_form}.
The original work considers the infinite time horizon impulse control problem and gives the analytical solution in an unbounded domain.
We first restrict the domain to be bounded, and add the boundary condition that gives a solution equivalent to that for the unbounded domain.
In Section \ref{sec:ofh_inf} we discuss the validity of this boundary condition using a numerical experiment.
Our main target problem, the finite time horizon optimal forest harvesting problem for which a solution has not been obtained analytically, is discussed in Section \ref{sec:ofh_fin}.

\subsection{Optimal forest harvesting problem}
\label{sec:ofh_form}
Let $X_t $ be the biomass of a forest at time $t$ and let $\tau_1<\tau_2<\cdots$ be the harvesting times.
We cut all trees in the forest at time $\tau_i$ and replant the biomass $\tilde{x} \in \mathbb{R}$.
Suppose that the growth of the biomass follows a geometric Brownian motion so that $X_t$ is governed by
\begin{align}
	&\begin{cases}
		\displaystyle dX^v_t = \mu X^v_t dt + \sigma X^v_t dW_t, \quad \tau_j < t < \tau_{j+1},\\
		X^v_{\tau_j} = \tilde{x}, \quad j=1,2,\cdots,\\
	\end{cases}
	 \label{eq:def_biomass_process}
\end{align}
where $\mu$ and $\sigma$ are positive constants.
Furthermore, we suppose that $\tau_i$ satisfies the conditions to be the intervention time, that is, $\tau_i$ is an $\mathcal{F}_t$-stopping time and $\tau_i<\infty$ almost surely. 
Then $X^v_t$ has a unique strong solution and $v:=(\tau_1,\tau_2,\cdots)$ is the admissible impulse control strategy.

The original problem formulated by Willassen \cite{Willassen1998} is defined as the infinite time horizon optimal impulse control problem.
Let $\beta\in(0,1)$ be the proportional harvesting cost and $Q>0$ be the replanting cost of the biomass $\tilde{x}$. 
Then the performance criterion is 
\begin{equation}
\label{eq:inf_perf_cri}
J^v(x) = \mathbb{E} \left[ \left.\sum_{i=1}^\infty e^{-\lambda \tau_i} \left( (1-\beta)X_{\tau_j -} - Q \right) \right| X^v_0=x  \right],
\end{equation}
where $\lambda>0$ is the discounting factor.
We impose on $\tilde{x}$, $\beta$ and $Q$ the condition that $(1-\beta)\tilde{x}<Q$:
 if this is not the case then the optimal strategy is that we harvest trees immediately after replanting which is a vacuous solution.
Harvesting immediately after replanting also suggests that $x^R=\tilde{x}$ for any $t \in (0,\infty)$.

The value function corresponding to the criterion (\ref{eq:inf_perf_cri}) is defined by $V(x)=\sup_{v}J^v(x)$, $x\in\mathbb{R}$,
and hence the corresponding HJBQVI is given by 
\begin{equation}
\max \left(
 \frac{\sigma^2 x^2}{2} \partial^2_x V(x)+ \mu x \partial_x V(x)- \lambda V(x) ,
 V(\tilde{x}) + (1-\beta)x - Q - V(x)
\right) = 0.
\label{eq:inf_qvi}
\end{equation}
Willassen solved the HJBQVI (\ref{eq:inf_qvi}) explicitly:
\begin{align}
	V(x) &=
	\begin{cases}
		\Psi(x) & \mathrm{for}\;x<y,\\
		(1-\beta)x- Q +\Psi(\tilde{x}) & \mathrm{for}\;x\geq y,
	\end{cases}
	 \label{eq:inf_ana_sol}
\end{align}
where 
\[
\Psi(x) =  \frac{(1-\beta)y}{\gamma}\left( \frac{x}{y} \right)^\gamma,
\quad
\gamma = \frac{\sigma^2-2\mu + \sqrt{(\sigma^2-2\mu)^2+8\sigma^2\lambda}}{2\sigma^2}
\]
and $y > \tilde{x}$ is a solution of
\[
y = \frac{\gamma Q - (1-\beta)y\left(\tilde{x}/y\right)^\gamma }{(1-\beta)(\gamma-1)}.
\]
We call $y$ the strategy switch point.
The key ideas for obtaining the solution are as follows:
(i) the condition for the cost suggests that we should wait to harvest until $X_t$ exceeds a certain value $y$;
(ii) because the partial differential equation $\frac{\sigma^2 x^2}{2} \partial^2_x V(x) + \mu x \partial_x V(x) - \lambda V(x)=0$
	has a general analytical solution, we obtain the value function by connecting this analytical solution and $V(\tilde{x}) + (1-\beta)x - Q $ smoothly at $y$.

We extend the above problem by considering the case that the farmer exits the forestry business at time $T$, meaning that he harvests the all trees and does not replant at time $T$.
Then the problem becomes a finite time horizon problem and the performance criterion is modified as follows:
\begin{equation}
\label{eq:fin_perf_cri}
J^v_t(x) = \mathbb{E} \left[ \left.\sum_{t<\tau_i<T} e^{-\lambda \tau_i} (1-\beta)X^v_{\tau_j^-} 
	+ e^{-\lambda T} \left( (1-\beta)X^v_{T} - Q \right) \right| X_t = x \right].
\end{equation}
The value function corresponding to this problem is defined by $V_t(x) = \sup_{v} J^v_t(x)$, $t\in[0,T)$, $x\in\mathbb{R}$,
and thus the corresponding HJBQVI is given by
\begin{align}
&\max \left(
	\partial_t V_t(x) + \frac{\sigma^2 x^2}{2} \partial^2_x V_t(x) + \mu x \partial_x V_t(x),
	\right. \nonumber \\ 
&\qquad\qquad\left.
	V_t(\tilde{x}) + e^{-\lambda t} \left( (1-\beta)x - Q \right) - V_t(x)
\right) = 0
\label{eq:fin_qvi}
\end{align}
with terminal condition $V_T(x)=e^{-\lambda T} (1-\beta)x$.

In this case an analytical solution is not available. However, we can anticipate the behavior of the solution of the HJBQVI (\ref{eq:fin_qvi}).
Because the performance criterion $J^v_0$ is equivalent to that for the infinite time horizon case (\ref{eq:inf_perf_cri}) in the limit as $T\rightarrow\infty$,
the value function $V_t$ and the optimal control $\hat{v}$ coincide with those for the infinite horizon provided $T$ is large enough and $t \ll T$.

\subsection{Determination of the bounded domain and boundary condition}
\label{sec:ofh_inf}
The idea of introducing the strategy switch point $y$ is significant for determining the candidate finite domain and boundary condition.
Let $x_{\max}$ be a positive real value, let $\mathcal{S}=(0,x_{\max})$ be the candidate bounded domain 
 and let the boundary $\partial \mathcal{S}=\{0,x_{\max}\}$ be the candidate boundary.
We first define $\psi$ as the function giving the boundary condition in the infinite time horizon case such that
\begin{align}
	\begin{cases}
		\psi(0) = 0,\\
		\psi(x_{\max}) = V(\tilde{x}) + (1-\beta)x_{\max}- Q.
	\end{cases}
	\label{eq:inf_boundary}
\end{align}
Because we cannot expect forest growth after the biomass reaches 0, the value function should be 0 at $x=0$.
We should expect that this boundary condition will give the same solution in the infinite domain case if $x_{\max}$ is sufficiently larger than $y$.

We examine this candidate domain and boundary condition by solving the corresponding HJBQVI (\ref{eq:inf_qvi}) numerically.
The method for solving the infinite horizon HJBQVI was proposed by {\O}ksendal and Sulem \cite{oksendal2007applied}.
Because our state process is a 1-dimensional geometric Brownian motion, we are able to use their method with no extra assumptions.
The HJBQVI is discretized using the standard finite difference scheme with a central difference,
 and we denote by $\delta_x$ be the finite difference step and by $N^x$ the number of grid points.
We determine the parameters from Table \ref{tb:parameters}.
In this situation, the value of the strategy switch point is $y=5.495503$.
\begin{table}[htbp]
	\centering
	\begin{tabular}{|l|l|l|}
		\hline
		Parameter  & Description           & Value \\ \hline
		$x_{\max}$ & right limit of the domain & 10    \\
		$\tilde{x}$ & initial & 1    \\
		$\beta$ & harvesting cost rate & 10\%    \\
		$Q$ & replanting cost & 2    \\
		$\mu$ & expected growth rate & 1    \\
		$\sigma$ & volatility & 1    \\
		$\lambda$ & discount factor & 2    \\
		\hline
	\end{tabular}
	\caption{Parameters}
	\label{tb:parameters}
\end{table}

The numerical results obtained using the algorithm are as follows.
\begin{figure}[htbp]
	\begin{minipage}[t]{0.48\columnwidth}
		\begin{center}
			\includegraphics[height=5.5cm]{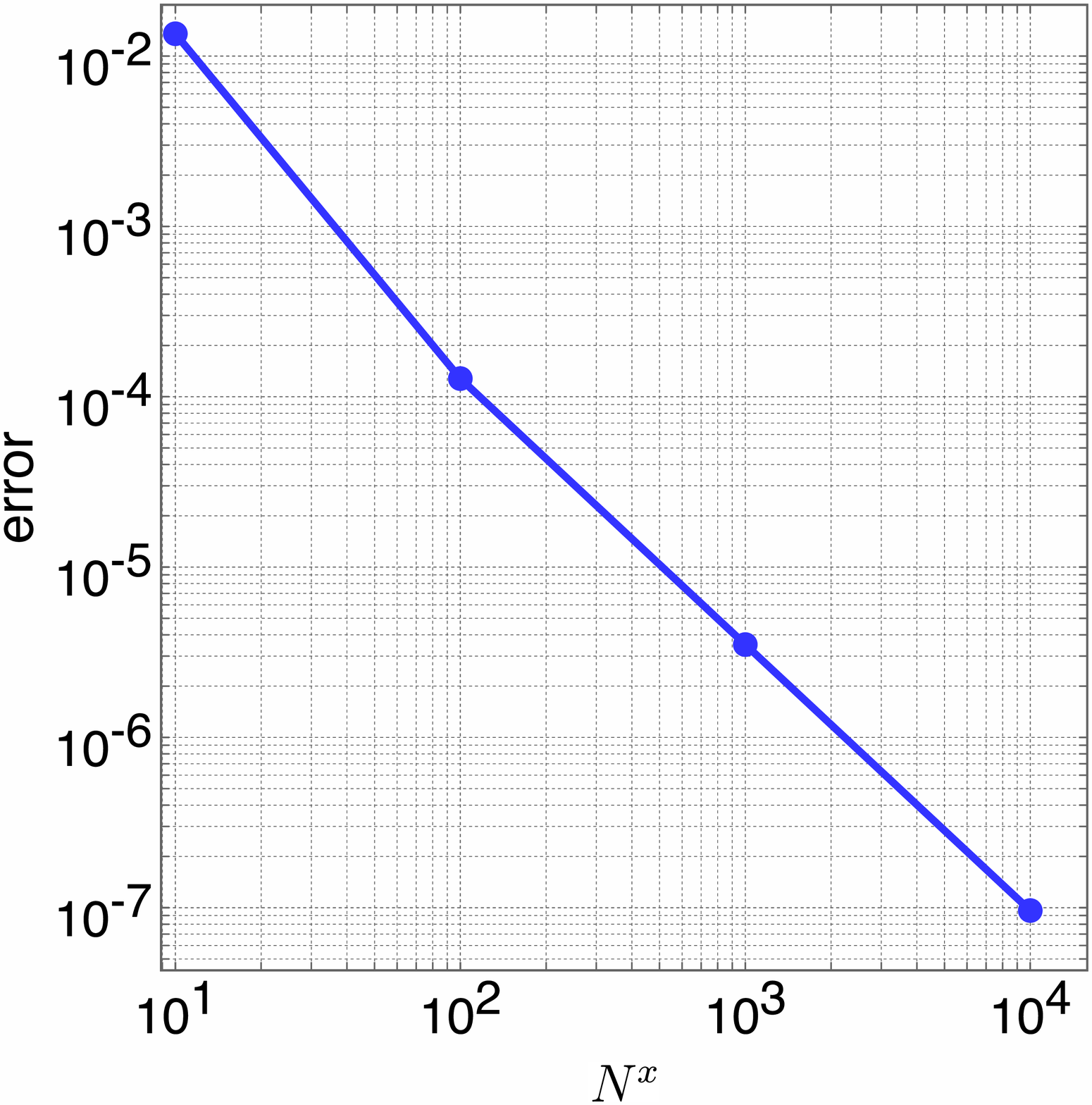}
			\caption{Maximum error. }
			\label{fig:inf_error}
		\end{center}
	\end{minipage}
	\hfill
	\begin{minipage}[t]{0.48\columnwidth}
		\begin{center}
			\includegraphics[height=5.5cm]{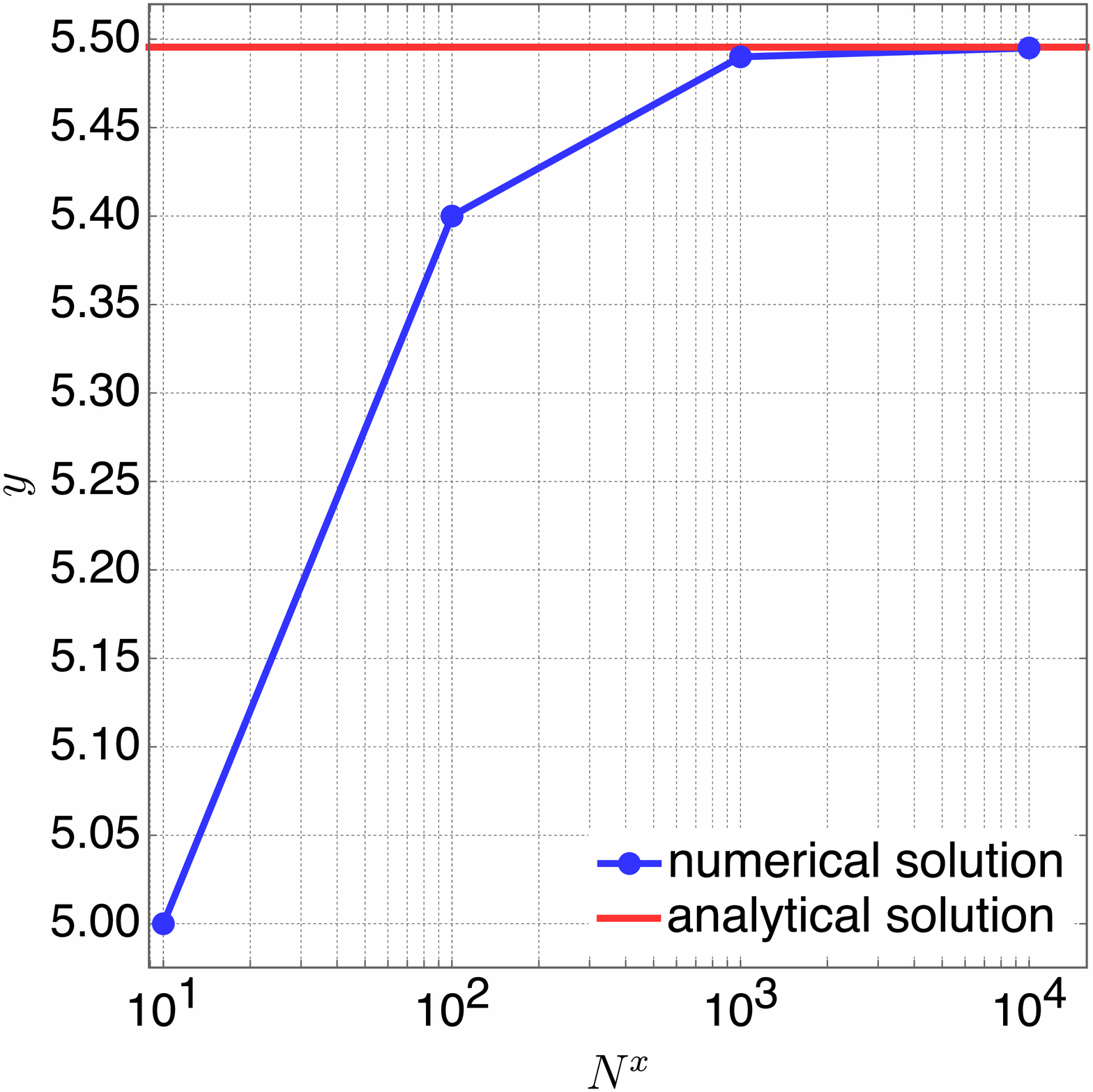}
			\caption{Strategy switch point. }
			\label{fig:inf_switch}
		\end{center}		
	\end{minipage}
\end{figure}
Figure \ref{fig:inf_error} shows the maximum error of the value function when comparing the analytical solution (\ref{eq:inf_ana_sol}) and the numerical solution.
The order of the error is $\mathcal{O}(\delta_x^2)$, which coincides with that implied by the finite difference scheme with a central difference.
Figure \ref{fig:inf_switch} displays the analytical strategy switch point with a red line and the numerical strategy with a blue line,
 indicating that the numerical results with the switch point agree with the $\delta_x$-order.
Therefore we conclude that our candidate domain and boundary condition are valid.

We next define the boundary condition in the finite time horizon case to be slightly modified from that in the infinite time horizon case:
\begin{align*}
	\begin{cases}
		\psi(t,0) = 0,\\
		\psi(t,x_{\max}) = V_t(\tilde{x}) + e^{-\lambda t} \left( (1-\alpha)x_{\max} -Q \right).
	\end{cases}	
\end{align*}
This boundary condition works well in the case that $x_{\max}$ is sufficiently large: the value of $x_{\max}$ should be larger than in the infinite time horizon case.
Because of the replanting cost $Q$, we expect that the optimal strategy close to the terminal time $T$ is to not cut down the trees.
However our boundary condition enforces cutting at the boundary point $x_{\max}$.
We can avoid this contradiction by taking the value of $x_{\max}$ to be large enough,
so that the contribution of the cost $Q$ to the value function $V_t(x)$ vanishes.
The reasons for this are that
 the harvest profit increases with the growth of the biomass $x$, and
 the cost $Q$ does not depend on the biomass $x$.
We discuss this issue again in the following section with reference to the numerical results.

At the end of this subsection, we mention the computational load for solving the infinite time horizon HJBQVI using this algorithm.
\begin{figure}[htbp]
	\begin{minipage}[t]{0.48\columnwidth}
		\begin{center}
			\includegraphics[height=5cm]{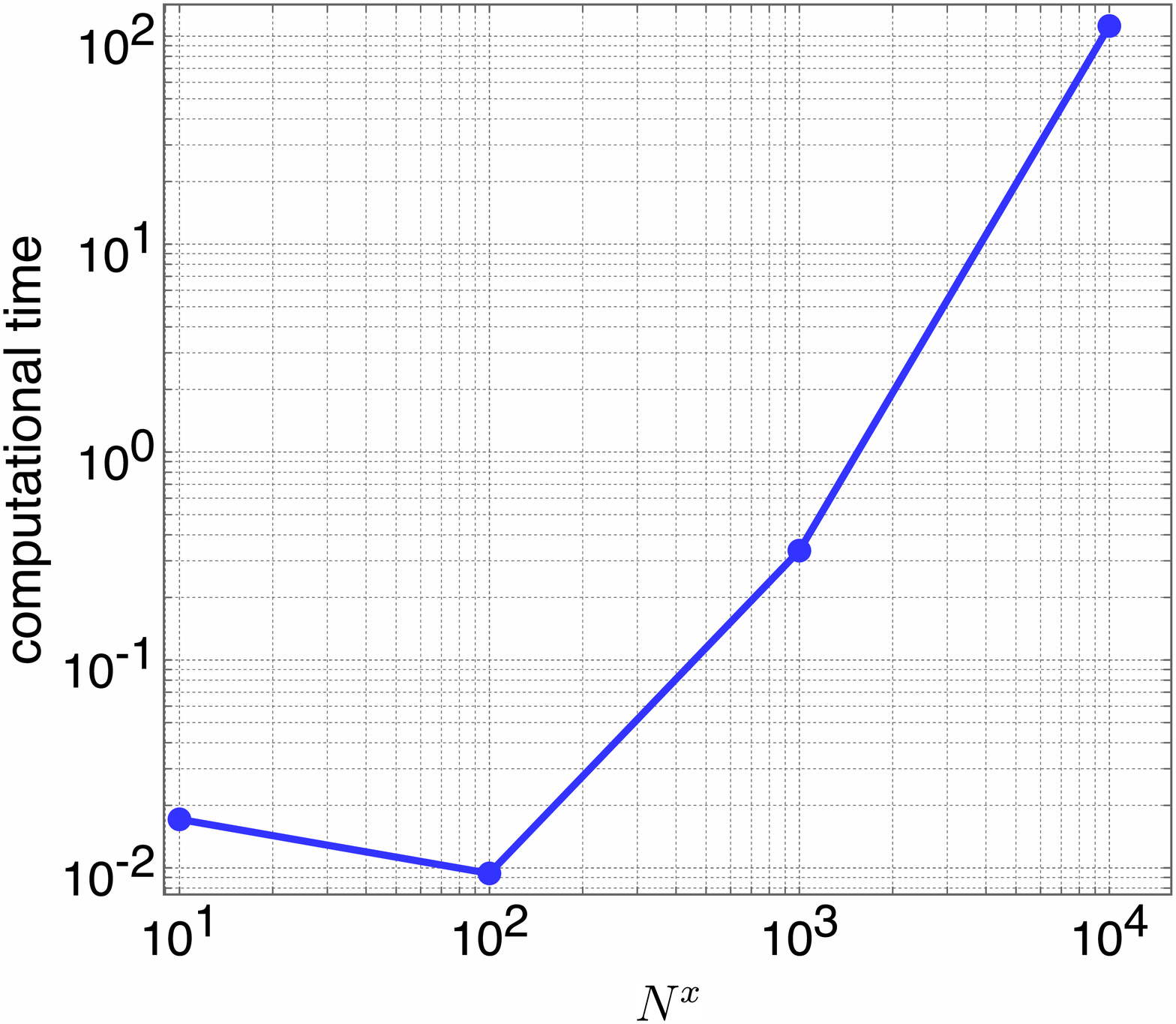}
			\caption{Computational time. }
			\label{fig:inf_time}
		\end{center}		
	\end{minipage}
	\hfill
	\begin{minipage}[t]{0.48\columnwidth}
		\begin{center}
			\includegraphics[height=5cm]{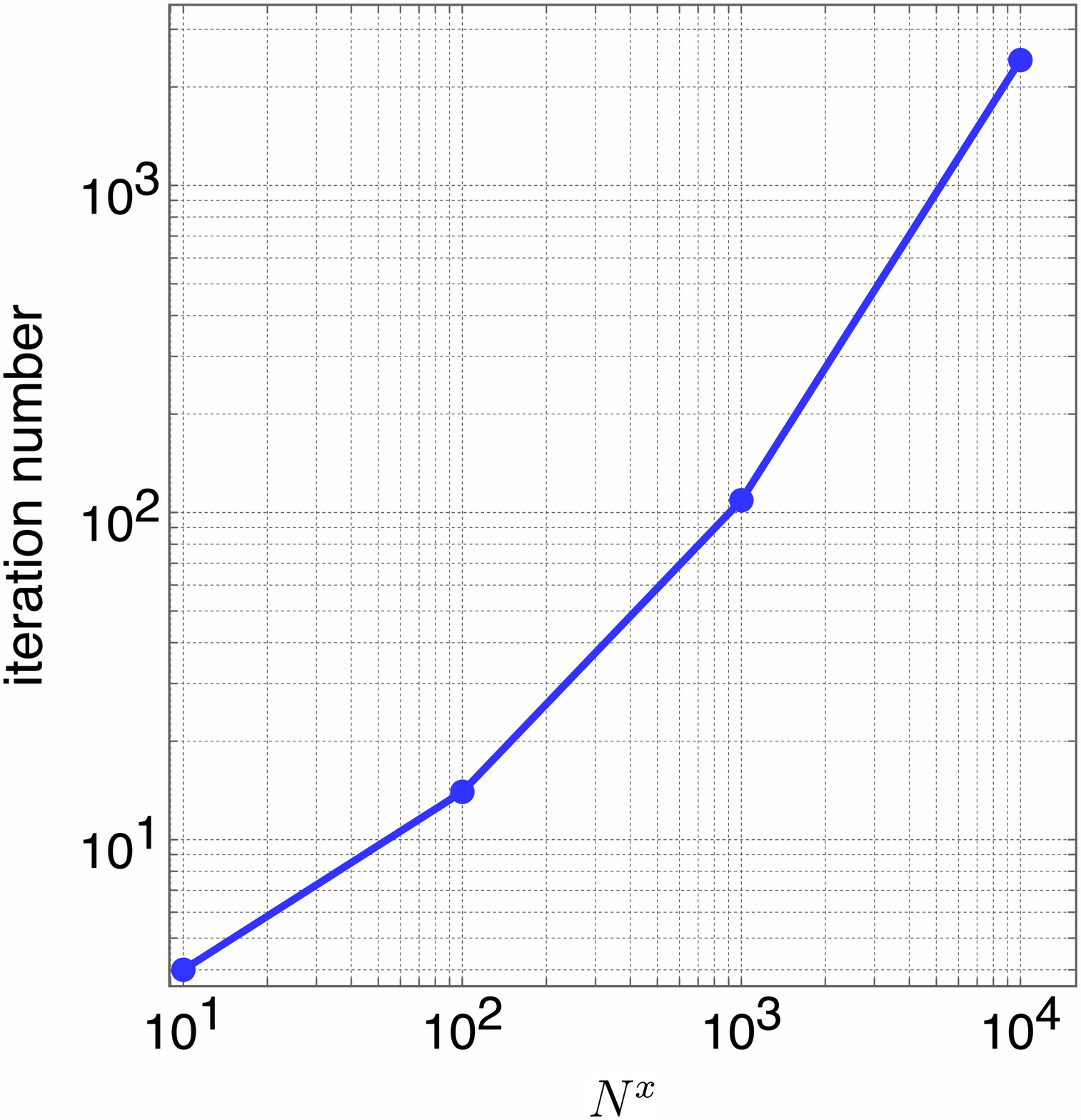}
			\caption{Number of iterations. }
			\label{fig:inf_itr}
		\end{center}
	\end{minipage}
\end{figure}
Figure \ref{fig:inf_time} describes the computational time as a function of the number of grid points
\footnote{
	The details of our computational resources are as follows.
	The computer we use has an Intel Core i5 650 @3.20GHz and 4GB RAM.
	Our code is parallelized using OpenMP, and we use PARDISO \cite{Schenk2004} for solving linear equations.
}.
We see that the computational time grows exponentially as $\delta_x$ becomes finer.
The main load is the growth in the size of the linear equations, which is inevitable until we use the finite difference scheme.
The search for the optimal regular and impulse controls on each grid point is the second load: the index $i$ takes a larger range of values when $N^x$ increases.
There is the possibility that we can reduce this factor through massive parallel computing such as GPGPU.
Figure \ref{fig:inf_itr} shows the inherent load of this algorithm.
Since the finer finite difference step allows us to compute a more accurate value function, we need more iterations to reach the fixed point.

\subsection{Finite time horizon case}
\label{sec:ofh_fin}
We now treat the main issue of this section.
In this situation, the matrix operators in our algorithm are defined as follows:
\begin{align*}
	&\displaystyle L^{a,k}_{ij}=\begin{cases}
		\displaystyle -\frac{\sigma^2 x_i^2}{\delta_x^2} &i=j,\\
		\displaystyle \frac{\sigma^2 x^2}{2\delta_x^2} \pm \frac{\mu x}{2\delta_x} &i=i\pm1, \\
		\displaystyle 0  &\mathrm{otherwise},
	\end{cases}\\
	&f^{a,k} = 0, \quad \eta(i,z_i) = j\in\mathbb{N}\;\mathrm{s.t.}\;x_j=\tilde{x},\\
	&K^{z,k}_i = e^{-\lambda t_k} \left( (1-\alpha)x_i - Q \right).
\end{align*}
We set $T=3.0$ and $N^t=3000$.
The other parameters have the same values as in Table \ref{tb:parameters}, except for $x_{\max}$.

\begin{figure}[htbp]
	\begin{center}
		\includegraphics[height=5cm]{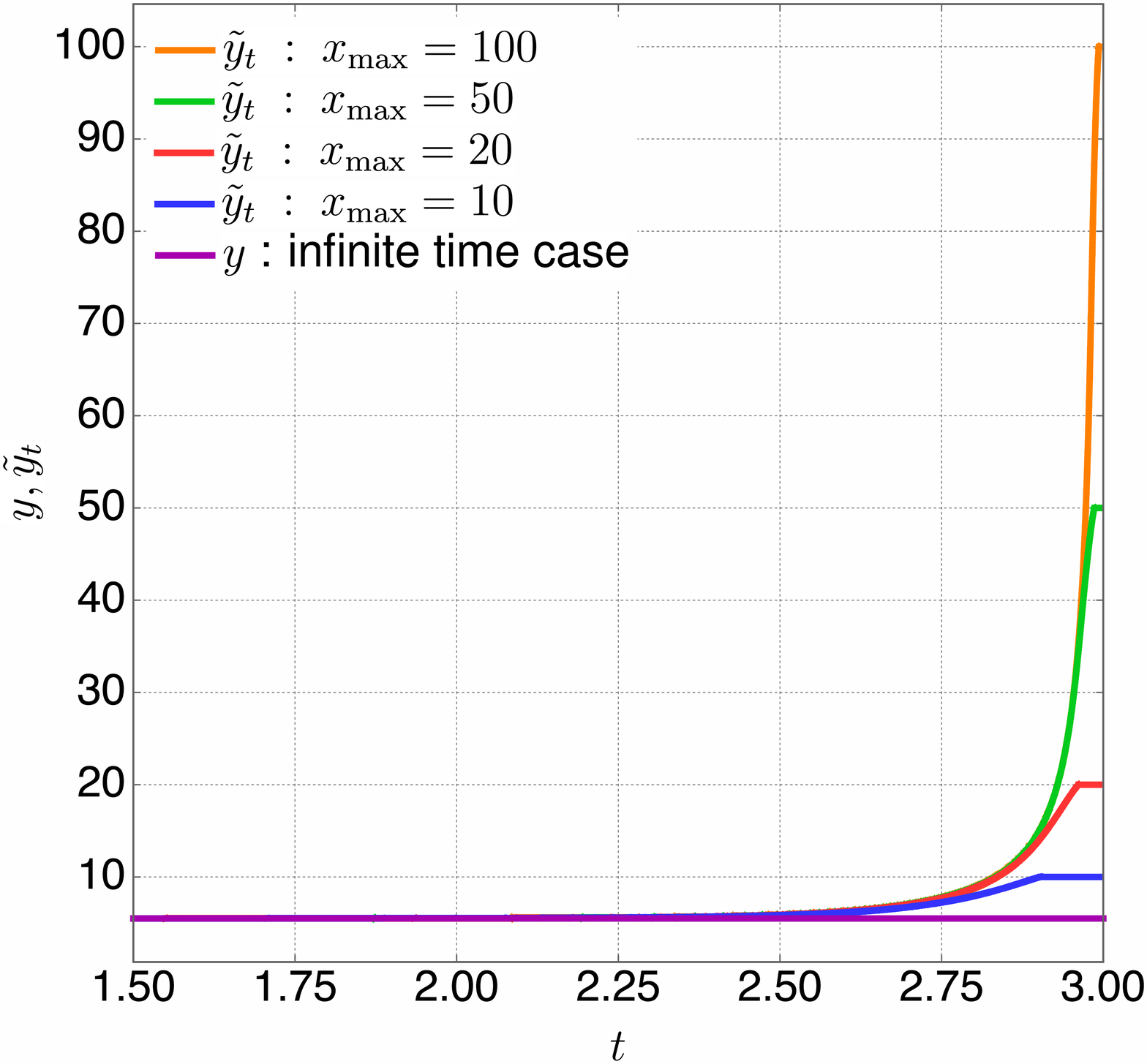}
		\caption{$\tilde{y}_t$ for various values of $x_{\max}$.
			 The orange, green, red and blue lines indicate the strategy switch point obtained with the boundaries $x_{\max}=$100, 50, 20 and 10 respectively.
			 The purple line indicates the infinite time horizon case.}
		\label{fig:fin_switch}
	\end{center}
\end{figure}

The first numerical result we discuss in this subsection is the behavior of $\tilde{y_t}$, the strategy switch point at time $t$.
As mentioned in Section \ref{sec:ofh_form}, we expect that $\tilde{y_t}$ should converge to $y$, the strategy switch point for the infinite time horizon case,
 if $t$ approaches $0$ and $T$ is large enough.
We also recall that $\tilde{y_t}$ contains error if $x_{\max}$ is not large enough.
Hence, we examine the various values of $x_{\max}$. The results are shown in Figure \ref{fig:fin_switch}.
We first mention that, remarkably, $\tilde{y}_t$ converges to $y$ regardless of the value of $x_{\max}$.
This coincides with our expectation discussed in the previous subsection.
We can understand the behavior of $\tilde{y}_t$ as follows.
The behavior of $\tilde{y}_t$ close to $T$, which takes a much larger value than $\tilde{y}$, suggests that 
 we should keep the trees.
If the time $t$ disengages from $T$, the benefit of the waiting to harvest caused by the replanting cost $Q$ vanishes and hence the value of $\tilde{y}_t$ decreases.

We next discuss the error due to the improper boundary condition.
In the previous subsection we hypothesized that the error vanishes. 
The behavior of $\tilde{y}_t$ gives evidence supporting this hypothesis as follows.
First, the length of the time interval for which $\tilde{y}_t$ is fixed at $x_{\max}$ near the terminal time $T$.
The length of this time interval decreases as $x_{\max}$ increases, and we cannot distinguish the interval in the case that $x_{\max}=100$ in Figure \ref{fig:fin_switch}.
Second, $\tilde{y}_t$ behave consistently when $t \ll T$.
When $t$ approaches $T$, $\tilde{y}_t$ begin to behave inconsistently with regard to the value of $x_{\max}$.

\begin{figure}[htbp]
	\begin{minipage}[t]{0.48\columnwidth}
		\begin{center}
			\includegraphics[height=5cm]{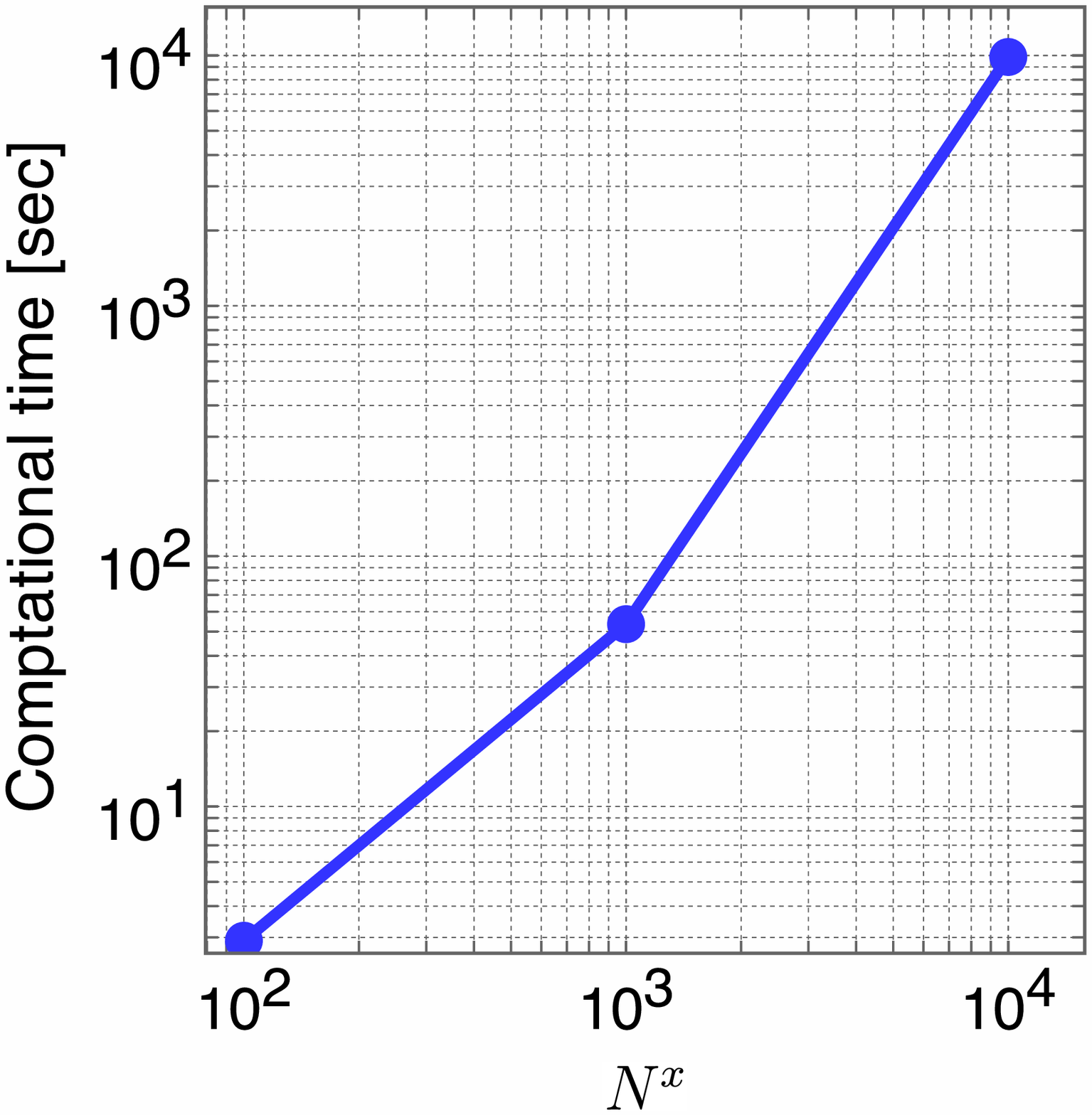}
			\caption{Computational time. }
			\label{fig:fin_time}
		\end{center}		
	\end{minipage}
	\hfill
	\begin{minipage}[t]{0.48\columnwidth}
		\begin{center}
			\includegraphics[height=5cm]{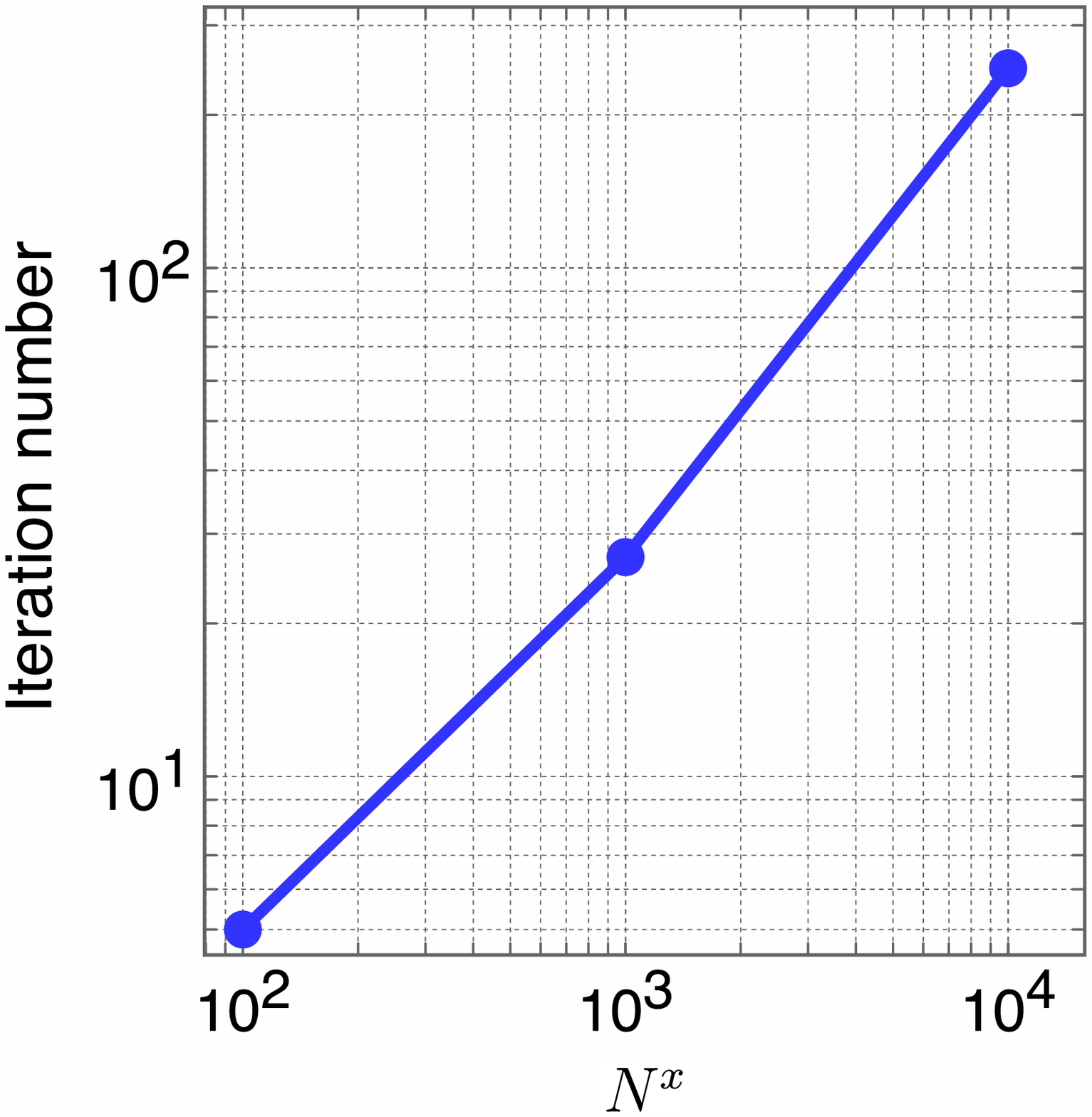}
			\caption{Number of iterations. }
			\label{fig:fin_itr}
		\end{center}
	\end{minipage}
\end{figure}

Finally, we discuss the computational load.
We set $x_{\max}=100$ and display the results in Figure \ref{fig:fin_time} and Figure \ref{fig:fin_itr}.
Figure \ref{fig:fin_time} describes the computational time where the horizontal axis is the number of grid points.
The new load factor is obviously the time grid size $N^t$.
The growth in the size of the linear equations and the search cost for the optimal controls $\hat{a}^k$ and $\hat{z}^k$ are the same  as in the previous subsection.
However, the number of iterations for convergence to the fixed point is slightly different.
Figure \ref{fig:fin_itr}, which shows that the maximum number of iterations for each time step, suggests that
 we only need approximately 10\% of the iterations compared with the same grid size in the infinite time horizon case.
This is due to the determination of the initial value of $\phi^k$.
Because $\delta_t$ is small enough, the difference between $\phi^k$ and $\phi^{k+1}$ is expected to be small.
Hence we are able to reduce the number of iterations from the infinite case where we have no information about the solution.

\section{Summary}
We have proposed a new numerical method for solving the HJBQVI associated
 with the combined impulse and stochastic optimal control problem over a finite time horizon.

The key points of our numerical scheme are that
 (i) we discretize the HJBQVI using the forward difference for the time grid;
 (ii) the impulse control part is kept in the same time step.
Hence the present method is an implicit method.
The main advantage of the present method is that it is independent of the specific model formulation.
The present work gives a proof of stability in the general framework and hence
 our method can be applied without a stability proof.
We have provided a detailed procedure for the algorithm in matrix form for the sake of easy implementation.

We also apply our method to the optimal forest harvesting problem.
The original problem formulated by Willassen \cite{Willassen1998} is defined as an infinite time horizon impulse control problem
 and has an analytical solution.
We introduce a terminal time representing the time of exit from the forestry business
 which turns the above problem into a finite time horizon problem.
Because the original problem is defined on an unbounded domain,
 we introduce an equivalent bounded domain and boundary condition 
 and verify them through a numerical experiment.

The analytical solution of our finite time horizon problem is unavailable and hence we solve it using our algorithm.
The behavior of the obtained optimal strategy is reasonable:
 the strategy coincides with the infinite time horizon strategy when the terminal time approaches infinity;
 the strategy switch point, the threshold for harvesting the trees, is much higher than in the infinite horizon case near the terminal time.

\appendix

\section{Stability}
\label{sec:stability}
We define $A \in \mathbb{R}^{N^x \times N^x}$ and $b, \varphi^k \in \mathbb{R}^n$ as
	\begin{align*}
		&A_{ij} = \begin{cases}
			\delta_t L^{a^*,k}_{ij} - \mathbf{1}_{ \{i=j\} } & \mathrm{if}\; i \in \mathcal{I}^{*^k},\\
			M^{z^*}_{ij} - \mathbf{1}_{ \{i=j\} } & \mathrm{if}\; i \in \{1,\cdots,N^x\}\backslash \mathcal{I},
		\end{cases}\\
		&b_i = \begin{cases}
			\displaystyle \delta_t f^{a^*,k}_i + \varphi^{k+1}_i
					 + \sum_{j=N^x+1}^{N^x+\bar{N}^x} \delta_t L^{a^*,k}_{ij} \psi(t_k,x_j)  & \mathrm{if}\; i \in \mathcal{I}^{*^k},\\
			K^{z^*,k}_i & \mathrm{if}\; i \in \{1,\cdots,N^x\}\backslash \mathcal{I},
		\end{cases}\\
		&\varphi^k_i = \phi^k_i, \qquad i \in \{1,\cdots,N^x\},
	\end{align*}
where
\begin{align*}
&a^*  = \argmax_a \left\{ \frac{\phi^{k+1} -\phi^k}{\delta_t} + L^{a,k} \phi^k  + f^{a,k}  \right\}, \\
&z^* =  \argmax_z \left\{ M^{z}\phi^k  +K^{z,k} -\phi^k \right\}, \\
&\mathcal{I}^{*k} = \left\{ i \in \{1,\cdots,N^x\} \left| 
		\frac{\phi^{k+1}_i -\phi^k_i}{\delta_t} + \left( L^{a^*,k}\phi^k \right)_i  + f^{a^*,k}_i 
			 > \left( M^{z^*}\phi^k \right)_i +K^{z^*,k}_i -\phi^k_i \right. \right\}.
\end{align*}
Then the $N^x$-dimensional linear equation $A\varphi+b=0$ is equivalent to the matrix QVI (\ref{eq:conv_matrix_qvi}).
Our goal is to show that $A$ is invertible and that $\Vert \varphi^k \Vert_\infty \leq C$ where $C$ is a constant independent of
 the grid sizes $\delta_t$ and $\delta$.

We first consider the case that $\mathcal{I}^{*k} = \{ 1,\cdots,N^x \}$.
By condition (\ref{eq:ass_central}), we find that $A_{ii}<0$, $A_{ij} \geq 0$ and
$\sum_{j} A_{ij} < -1$ for $i,j=1,\cdots,N^x$ and $j \neq i$. 
This implies that $A$ is an M-matrix and hence $A$ is invertible.
Next, we define $\tilde{i} \in \{1,\cdots,N^x\}$ as $\tilde{i} = \argmax_i \left| \varphi^k_i \right| $.
Because $A_{ij} \geq 0$ for $i \neq j$
 and $\Vert \varphi^k \Vert_\infty \geq \mathrm{sign}(\varphi^k_{i}) \varphi^k_{j}$ for $i,j\in \{1,\cdots,N^x \}$,
 we have
\begin{align*}
\sum_{j \neq \tilde{i}} A_{\tilde{i} j} \Vert \varphi^k \Vert_\infty  \geq & \sum_{j \neq \tilde{i}} A_{\tilde{i} j} \mathrm{sign}(\varphi^k_{\tilde{i}}) \varphi^k_j \\
\Leftrightarrow 
- \sum_{j \neq \tilde{i}} A_{\tilde{i} j} \Vert \varphi^k \Vert_\infty  \leq & - \mathrm{sign}(\varphi^k_{\tilde{i}}) \sum_{j \neq \tilde{i}} A_{\tilde{i} j}  \varphi^k_j \\
&= \mathrm{sign}(\varphi^k_{\tilde{i}}) \left( A_{\tilde{i} \tilde{i}} \varphi^k_{\tilde{i}} + b_{\tilde{i}} \right) \\
\Leftrightarrow - \sum_{j \neq \tilde{i}} A_{\tilde{i} j} \Vert \varphi^k \Vert_\infty -&  \mathrm{sign}(\varphi^k_{\tilde{i}}) A_{\tilde{i} \tilde{i}} \varphi^k_{\tilde{i}}
		  \leq \mathrm{sign}(\varphi^k_{\tilde{i}}) b_{\tilde{i}}.
\end{align*}
Evaluating the left- and right-hand sides of the last inequality:
\begin{align*}
- \sum_{j \neq i} A_{\tilde{i} j} \Vert \varphi^k \Vert_\infty -  \mathrm{sign}(\varphi^k_{\tilde{i}}) A_{\tilde{i} \tilde{i}} \varphi^k_{\tilde{i}}  
&= - \sum_{j} A_{\tilde{i} j} \Vert \varphi^k \Vert_\infty  +  A_{\tilde{i} \tilde{i}} \Vert \varphi^k \Vert_\infty
	 -  \mathrm{sign}(\varphi^k_{\tilde{i}}) A_{\tilde{i} \tilde{i}} \varphi^k_{\tilde{i}} \\ 
&= - \sum_{j} A_{\tilde{i} j} \Vert \varphi^k \Vert_\infty  +  A_{\tilde{i} \tilde{i}} ( \Vert \varphi^k \Vert_\infty -  \mathrm{sign}(\varphi^k_{\tilde{i}}) \varphi^k_{\tilde{i}} ) \\
&= - \sum_{j} A_{\tilde{i} j} \Vert \varphi^k \Vert_\infty \\
&\geq  \Vert \varphi^k \Vert_\infty;
\end{align*}
\[
\mathrm{sign}(\varphi^k_{\tilde{i}}) b_{\tilde{i}} \leq |b_{\tilde{i}}|
 \leq \left| \delta_t f^{a^*,k}_{\tilde{i}} \right| + \left| \varphi^{k+1}_{\tilde{i}} \right|
					 + \sum_{j=N^x+1}^{N^x+\bar{N}^x} \left| \delta_t L^{a^*,k}_{{\tilde{i}}j} \psi(t_k,x_j) \right|.
\]
Hence we obtain
\[
\Vert \varphi^k \Vert_\infty \leq \left| \delta_t f^{a^*,k}_{\tilde{i}} \right| + \left| \varphi^{k+1}_{\tilde{i}} \right|
					 + \sum_{j=N^x+1}^{N^x+\bar{N}^x} \left| \delta_t L^{a^*,k}_{{\tilde{i}}j} \psi(t_k,x_j) \right|.
\]
From the definition of the constants and functions, we find that $\Vert \varphi^k \Vert_\infty \leq C$,
 where C is a constant independent of $\delta_t$ and $\delta$.

We next consider the general $\mathcal{I}^{*k}$.
Let $\eta^{l,z,k}(i)$ be the $l$-th decomposition of the function $\eta$ with control set $z$ at time $t_k$ and let
$\bar{l}^{z,k}(i) = \inf \left\{ l \geq 1 \left| \eta^{l,z,k}(i) \in \mathcal{I}^{*k}  \right. \right\}$.
By equation (\ref{eq:ass_gamma}), we find that 
\begin{equation}
\varphi^k_i = \varphi^l_{\eta^{\bar{l}^{z^*,k}(i)}(i)} + \sum_{j=0}^{\bar{l}^{z^*,k}(i)-1}K^{z^*,k}_{\eta^{j,z^*,k}(i)} \; , \quad i \in \{1,\cdots,N^x \} \backslash \mathcal{I}^{*k}.
\label{eq:leq_imp_part}
\end{equation}
Hence the $N^x$-dimensional equation $A \varphi^k + b=0$ is reduced to the ($N^x-\#\mathcal{I}^{*k}$)- dimensional equation $\tilde{A}\varphi^k+\tilde{b}=0$,
 where $\#\mathcal{I}^{*k}$ is the number of elements in $\mathcal{I}^{*k}$.
Equation (\ref{eq:leq_imp_part}) implies that $\tilde{A}_{ii}<0$, $\tilde{A}_{ij} \geq 0$, and $\sum_{j} \tilde{A}_{ij} < -1$.
We find that $\tilde{A}$ is invertible and that $\Vert \varphi^k \Vert_\infty \leq C$ in the same manner as when $\mathcal{I}^{*k} = \{ 1,\cdots,N^x \}$.

\section{From the matrix QVI (\ref{eq:matrix_qvi}) to the fixed-point problem (\ref{eq:conv_matrix_qvi})}
\label{sec:conv_mat_fix}
We first verify that the matrix QVI (\ref{eq:matrix_qvi}) and the fixed-point problem (\ref{eq:conv_matrix_qvi}) are equivalent.
The equations (\ref{eq:bar_lf}) imply that
\begin{align*}
 &f^{a,k}_i		= \left( \frac{1}{\delta_t}+\frac{1}{h} \right) \bar{f}^{a,k}_i  - \frac{1}{\delta_t}\phi^{k+1}_i,\\
 &L^{a,k}		= \left( \frac{1}{\delta_t} + \frac{1}{h} \right)\bar{L}^{a,k} - \frac{1}{h} I.
\end{align*}
Hence the matrix QVI (\ref{eq:matrix_qvi}) can be rewritten in the following form:
\begin{align*}
&\begin{cases}
\max\left( 
	\displaystyle \left( \frac{1}{\delta_t} + \frac{1}{h} \right) \sup_{a\in\mathbb{U}^{N^x}} \left\{ 
		  \left( \bar{L}^{a,k}\phi^k \right)_i  + \bar{f}^{a,k}_i - \phi^k_i   \right\} \right.\\
	\displaystyle \left. \hspace*{2cm} ,\sup_{z\in\mathcal{Z}^{N^x}} \left\{ \left( M^{z}\phi^k \right)_i +K^{z,k}_i -\phi^k_i \right\}
\right)=0, \quad 1 \leq i \leq N^x,\\
\phi^k_i = \psi(t_k,x_i), \quad N^x < i \leq N^x+\bar{N}^x.
\end{cases}
\end{align*}
The parameters $\delta_t$ and $h$ are positive numbers, and thus this is equivalent to
\footnote{
	Note that $\max[ cf(x), g(x) ]=0$ is equivalent to $\max[ f(x), g(x) ]=0$
	for every $c>0$.
}
\begin{align*}
&\begin{cases}
\max\left( 
	\displaystyle\sup_{a\in\mathbb{U}^{N^x}} \left\{ 
		 \left( \bar{L}^{a,k}\phi^k \right)_i  + \bar{f}^{a,k}_i  \right\} - \phi^k_i \right.\\
	\displaystyle \left. \hspace*{2cm} ,\sup_{z\in\mathcal{Z}^{N^x}} \left\{ \left( M^{z}\phi^k \right)_i +K^{z,k}_i  \right\} -\phi^k_i
\right)=0, \quad 1 \leq i \leq N^x,\\
\phi^k_i = \psi(t_k,x_i;\phi), \quad N^x < i \leq N^x+\bar{N}^x.
\end{cases}
\end{align*}
Therefore we obtain the fixed-point problem (\ref{eq:conv_matrix_qvi}).

We next show that $\bar{L}^{a,k}$ is a contraction map displayed in matrix form.
By the definition of $h$, condition (\ref{eq:ass_central}) and equation (\ref{eq:bar_lf}), we find that
\begin{align}
&0\leq \bar{L}^{a,k}_{ij},<1  \quad i\in{1,\cdots,N^x}, \; j\in{1,\cdots,N^x+\bar{N}^x},\\
&\sum_{j=1}^{N^x+\bar{N}^x} \bar{L}^{a,k}_{ij} < 1, \quad  i\in{1,\cdots,N^x}.
\end{align}
Thus we obtain 
\begin{align*}
\left\Vert \bar{L}^{a,k}\phi^\prime - \bar{L}^{a,k}\phi \right\Vert_\infty
 &= \max_{1\leq i  \leq N^x} \left| \sum_{j=1}^{N^x+\bar{N}^x} \bar{L}^{a,k}_{ij}\phi^\prime_j - \bar{L}^{a,k}_{ij}\phi_j    \right| \\
 &\leq \max_{1\leq i  \leq N^x}  \sum_{j=1}^{N^x+\bar{N}^x}  \left| \bar{L}^{a,k}_{ij}\phi^\prime_j - \bar{L}^{a,k}_{ij}\phi_j    \right| \\
 &\leq \max_{1\leq j^\prime \leq N^x} \left| \phi^\prime_{j^\prime} - \phi_{j^\prime} \right| 
 		\max_{1\leq i  \leq N^x}  \sum_{j=1}^{N^x+\bar{N}^x}  \left| \bar{L}^{a,k}_{ij} \right| \\
 &< \max_{1\leq j^\prime  \leq N^x} \left| \phi^\prime_{j^\prime} - \phi_{j^\prime} \right| \\
 &= \left\Vert \phi^\prime - \phi\right\Vert_\infty.
\end{align*}
Hence $\bar{L}^{a,k}$ is a contraction map displayed in matrix form.

Finally, we show that $\bar{L}^{a,k}$ satisfies the discrete maximum principle, that is, that
$$
\bar{L}^{a,k}\phi^\prime - \bar{L}^{a,k}\phi \leq \phi^\prime - \phi \Rightarrow \phi^\prime - \phi \geq 0.
$$
To do so, we obtain a contradiction to the above statement.
Assume that $\phi^\prime_n - \phi_n < 0$ and that $\phi^\prime_i - \phi_i \geq 0$ for $i\neq k$. 
Then we have
\begin{align*}
(\bar{L}^{a,k}\phi^\prime)_n - (\bar{L}^{a,k}\phi)_n &= \sum_j \bar{L}^{a,k}_{nj}(\phi^\prime_j-\phi_j)\\
&= \sum_{j\neq n} \bar{L}^{a,k}_{nj}(\phi^\prime_j-\phi_j) + \bar{L}^{a,k}_{nn}(\phi^\prime_n-\phi_n)\\
&\geq \bar{L}^{a,k}_{nn}(\phi^\prime_n-\phi_n)\\
&> (\phi^\prime_n-\phi_n).
\end{align*}
Hence we obtain
$$
^\exists n \; \mathrm{s.t.} \; \phi^\prime_n-\phi_n <0 \; \Rightarrow \bar{L}^{a,k}\phi^\prime - \bar{L}^{a,k}\phi \nleq \phi^\prime - \phi,
$$
which establishes the statement.

\bibliography{/Users/mieda/Dropbox/bibtex/library,/Users/mieda/Dropbox/bibtex/text}

\end{document}